\documentclass[journal]{IEEEtran}
\usepackage{hyperref}
\hypersetup{
	colorlinks=true,
	linkcolor=black,     
	urlcolor=black,
	citecolor = black
}
\usepackage{bbm}
\usepackage{caption}
\usepackage{amsthm}
\usepackage{utfsym}
\usepackage{amssymb,graphicx,lineno}
\usepackage{rawfonts,amsfonts,amssymb,psfrag,color}
\usepackage{amsfonts}
\usepackage{amsmath,epsfig}
\usepackage{rawfonts,graphicx,amsfonts,amssymb}
\usepackage{bm}
\usepackage{algorithm, algorithmic}
\usepackage{subfig}
\usepackage{graphicx}
\usepackage{float}
\usepackage{verbatim}
\usepackage{color,xcolor}

\newtheorem{mypro}{Proposition}

\newtheorem{theorem}{Theorem}

\modulolinenumbers[5]
\usepackage{multirow}
\usepackage{booktabs}
\usepackage{cite}
\usepackage{color}

\usepackage{comment}

\begin{document}
	\title{A framework to generate sparsity-inducing regularizers for enhanced low-rank matrix completion}
	%
	\author{Zhi-Yong Wang and Hing Cheung So,~\IEEEmembership{Fellow,~IEEE}\thanks{E-mail: z.y.wang@my.cityu.edu.hk, hcso@ee.cityu.edu.hk. This work is supported by a grant from the Research Grants Council of the Hong Kong Special Administrative Region, China [Project No. CityU 11207922]. }}

	\maketitle
	\begin{abstract}
		Applying half-quadratic optimization to loss functions can yield the corresponding regularizers, while these regularizers are usually not sparsity-inducing regularizers (SIRs). To solve this problem, we devise a framework to generate an SIR with closed-form proximity operator.  Besides, we specify our framework using several commonly-used loss functions, and produce the corresponding SIRs, which are then adopted as nonconvex rank surrogates for low-rank matrix completion. Furthermore, algorithms based on the alternating direction method of multipliers are developed. 
		Extensive numerical results show the effectiveness of our methods in terms of recovery performance and runtime. 
		
	\end{abstract}
	\begin{IEEEkeywords}
		Matrix completion, rank minimization, rank surrogate, sparsity-inducing regularizer.
	\end{IEEEkeywords}
	\section{Introduction}
	\label{sec:intro}
	
	{L}{ow-rank} matrix completion (LRMC) aims to fill the unobserved entries of an incomplete matrix with the use of the low-rank property~\cite{Davenport2016}. It has numerous real-world applications such as image inpainting~\cite{WangZYATCSVT2023}, hyperspectral image restoration~\cite{WangZYRTCSVT2023} and collaborative filtering~\cite{ZengWJO2018}. That is because although these data have high-dimensional structure, their main information lies in a subspace with a much lower dimensionality. 	
	Roughly speaking, two strategies are widely used for LRMC, namely, matrix factorization~\cite{ShenY2014,ZhouTT2013} and rank minimization~\cite{Cand2009,Cai2010}. The former formulates the estimated matrix as a product of two much smaller matrices, but it requires knowing the prior rank information, which may be not easy to determine in real-world scenarios. 
	
	Different from matrix factorization, the rank minimization can estimate the rank of the observed matrix via employing an SIR to make the singular values sparse~\cite{LuC2016,LuC2014,LuC2015}. 
	One popular and feasible method is nuclear norm minimization (NNM)~\cite{Cand2009}. However, since NNM utilizes the $\ell_1$-norm to shrink all nonzero singular values by the same constant, the resultant solution is biased. 
	To solve this issue, nonconvex sparsity-inducing regularizers (SIRs) are suggested because they have less estimation bias than the $\ell_1$-norm~\cite{NieF2019}. 
	Many nonconvex SIRs are adopted to replace the $\ell_1$-norm for LRMC~\cite{LuC2016}, while they cannot ensure that the resultant subproblems associated with these SIRs are convex, and the closed-from solutions to these subproblems are not obtained, implying high computational costs. To handle this problem, a parameterized nonconvex SIR is proposed in~\cite{ParekhA2016}. However, to maintain the convexity of the subproblem associated with this SIR, one parameter is required to be properly set. Moreover, Gu $et$ $al$.~\cite{GuSH2017,GuSH2014} employ a weighted NNM (WNNM) for LRMC, yielding better low-rank recovery than the NNM.
	
	On the other hand, it has been analyzed that applying half-quadratic (HQ) optimization~\cite{NikolovaM2005} or the Legendre-Fenchel (LF) transform~\cite{RockafellarRT2004} to loss functions, including the Welsch, German-McClure (GMC) and Cauchy functions~\cite{MandanasFD2017}, can yield regularizers with closed-form proximity operators. However, these regularizers are usually not SIRs~\cite{WangZYRTSP2023}. In this work, our attempt is to answer an interesting and important question: \textit{Under what conditions, the resultant regularizers are SIRs} ?.
	
	Motivated by the results in~\cite{WangZYRTCSVT2023,WangZYRTSP2023,HeR2014,ChartrandR2012}, we devise a framework to generate SIRs with closed-form proximity operators. Note that the loss functions considered in~\cite{WangZYRTCSVT2023,WangZYRTSP2023} are a special case of our framework. Besides, it is analyzed that the subproblems associated with our SIRs are convex with closed-form solutions and our SIRs can yield a less bias solution than the $\ell_1$-norm.
	We then employ the SIRs for LRMC, and algorithms based on the alternating direction method of multipliers (ADMM) are developed.

	\section{Preliminaries}
	\label{sec:pro-for}
	\subsection{Proximity Operator}
	The Moreau envelope of a regularizer $\varphi(\cdot)$ is defined as~\cite{CombettesP2005}:
	\begin{equation}\label{Def_Pro}
		\begin{split}
			\min\limits_{y}~\frac{1}{2}(x-y)^2 + \lambda\varphi(y)
		\end{split}    	   	
	\end{equation}
	whose solution is given by the proximity operator:
	\begin{equation}\label{R-LSp}
		\begin{split}
			P_\varphi(x) := {\rm \arg}\min\limits_{y}~\frac{1}{2}(x-y)^2 + \lambda\varphi(y)
		\end{split}    	   	
	\end{equation} 
	In particular, when $\varphi(\cdot)=|\cdot|_1$, the solution to (\ref{Def_Pro}) is:
	\begin{equation}\label{Pro-L1}
		\begin{split}
			P_{\ell_1,\lambda}(x)={\rm max}\{0,|x|-\lambda\}{\rm sign}(x)
		\end{split}    	   	
	\end{equation}
	which is called the proximity operator of $|\cdot|_1$, also known as the soft-thresholding operator. 
	From (\ref{Pro-L1}), it is clear that the $\ell_1$-norm is an SIR. Here, $\varphi(\cdot)$ is called an SIR if the solution to (\ref{Def_Pro}) is zero if the magnitude of $|x|$ is no bigger than a threshold.
	On the other hand, although applying the HQ optimization to the Welsch, Cauchy or GMC function can give the corresponding regularizer with closed-form proximity operator~\cite{MandanasFD2017}, the generated regularizers are not SIRs.

	\subsection{Related Works}
	Given an incomplete matrix $\pmb X_\Omega\in \mathbb{R}^{m\times n}$ with ${\Omega} \subset \{1,\cdots,m\}\times\{1,\cdots,n\}$ being the index set of the observed elements, defined as:
	$$ \left[\pmb X_{\Omega}\right]_{ij}  = \left\{
	\begin{aligned}
		& X_{ij},\quad   {\rm if}~(i,j)\in \Omega  \\
		&0,\quad \quad   {\rm if}~(i,j)\in \Omega^c  \\
	\end{aligned}
	\right.
	$$
	where $\Omega^c$ is the complement of $\Omega$,
	LRMC can be solved by NNM~\cite{Cand2009}:
	\begin{equation}\label{nuclear-norm}
		\begin{split}
			&\mathop {\min}\limits_{\pmb M}~ \|{\pmb M}\|_*, ~ \text{s.t.} ~ \pmb M_{\Omega} = \pmb X_{\Omega}
		\end{split}    	   	
	\end{equation}
	where $\|\pmb M\|_*=\sum_{i=1}^{r} \pmb \sigma_i$ denotes the nuclear norm of the estimated matrix $\pmb M$ and $\pmb \sigma_i$ is the $i$th singular value of $\pmb M$. 
	Nevertheless, nuclear norm is equivalent to applying the $\ell_1$-norm to the singular value of a matrix, which shrinks all singular values with the same constant. However, as large singular values may dominate the main structure of real-world data, it is better to shrink them less~\cite{GuSH2017}. Hence nonconvex SIRs have been suggested to replace the $\ell_1$-norm as nonconvex rank surrogates for LRMC~\cite{LuC2016}:
	\begin{equation}\label{nonconvex-norm}
		\begin{split}
			&\mathop {\min}\limits_{\pmb M}~ \|{\pmb M}\|_\varphi, ~ \text{s.t.} ~ \pmb M_{\Omega} = \pmb X_{\Omega}
		\end{split}    	   	
	\end{equation}
	where $\|\pmb M\|_\varphi=\sum_{i=1}^{r} \varphi{(\pmb \sigma_i)}$ and $\varphi(\cdot)$ is a nonconvex SIR. As shown in~\cite{LuC2014,LuC2015}, solving (\ref{nonconvex-norm}) involves the proximity operator of $\varphi(\cdot)$. However, existing SIRs, like the $\ell_p$-norm with $0<p<1$ except for $p=\{1/2,2/3\}$~\cite{ShangFB2018}, may not have the closed-form proximity operator, implying that iterations are needed for its computation. 
	Furthermore, it is worth noting that when $\varphi{(\pmb \sigma_i)}=\pmb w_i \pmb \sigma_i$ with $\pmb w_i$ being a weight parameter, $\|\pmb M\|_\varphi$ is the weighted nuclear norm and (\ref{nonconvex-norm}) becomes the WNNM problem~\cite{GuSH2017}.

	\section{Framework to Generate SIR and its Application to LRMC}
	\label{sec:BVP}
	
	As the Huber~\cite{HeR2014}, truncated-quadratic~\cite{WangZYRTCSVT2023} and hybrid ordinary-Welsch (HOW)~\cite{WangZYRTSP2023} functions can yield SIRs using the LF transform, we generalize this class of functions and devise a framework to produce an SIR with closed-from proximity operator. Then, the generated SIR is considered as a nonconvex rank surrogate for LRMC.

	\subsection{Framework to Generate SIRs}
	\label{R-OR1MP}
	We first provide our framework in the following theorem.
	
	\begin{theorem}
		Consider a continuously differentiable loss function $\phi_{h,\lambda}(x)$, defined as:
		\begin{equation}\label{Def_framework}
			\phi_{h,\lambda}(x) = 
			\begin{cases}
				x^2/2, &|x|\leq \lambda\\
				a\cdot h(|x|)+b, &|x|\textgreater \lambda
			\end{cases}
		\end{equation}
		 where $a$ and $b$ are constants to make $\phi_{h,\lambda}(x)$ continuously differentiable, such that $g(x)=x^2/2-\phi_{h,\lambda}(x)$ is a proper, lower semi-continuous and convex function. Then it can be used to generate an SIR $\varphi_{h,\lambda}(\cdot)$ via LF transform, namely,
		\begin{equation}\label{covex_newphi}
			\phi_{h,\lambda}(x)=\mathop {\min}\limits_{y} ~ \frac{1}{2}(x-y)^2+ \lambda\varphi_{h,\lambda}(y)
		\end{equation}
		where $\varphi_{h,\lambda}(y)=\mathop {\max}\limits_{x} ~ \phi_{h,\lambda}(x)/\lambda-\frac{1}{2\lambda}(x-y)^2$.
		The solution to $y$ in (\ref{covex_newphi}) is given by the proximity operator:
		\begin{equation}\label{Pro_Ope}
			P_{\varphi_{h,\lambda}}(x)=\nabla g= {\rm max}\left\{0, |x|-a\cdot h'(|x|) \right\}\cdot {\rm sign}(x)
		\end{equation}
	\end{theorem}
	
	\emph{Proof:}
	The process of obtaining (\ref{Pro_Ope}) and (\ref{covex_newphi}) from (\ref{Def_framework})  is the same as that getting (16) and (17) from (9) in our previous work~\cite{WangZYRTSP2023}, and we omit the process due to page limit. 
	$\hfill\blacksquare$
	
	In addition, the expression of $\varphi_{h,\lambda}(\cdot)$ is generally unknown~\cite{NikolovaM2005,HeR2014}, and its properties are provided in the following proposition.
	\begin{mypro}\label{solution-proximal} 
		$\varphi_{h,\lambda}(\cdot)$ has the following properties:
		\begin{itemize}
			\item[(i)] Problem (\ref{covex_newphi}) is a convex problem although $\varphi_{h,\lambda}(y)$ is nonconvex.
			\item[(ii)] $P_{\varphi_{h,\lambda}}(x)$ is monotonically non-decreasing, namely, for any $x_1<x_2$, $P_{\varphi_{h,\lambda}}(x_1)\leq P_{\varphi_{h,\lambda}}(x_2)$.
			\item[(iii)] If $\phi_{h,\lambda}(x)$ is strictly concave for $x>\lambda$, the resultant proximity operator makes the solution have less bias than the proximity operator of the $\ell_1$-norm.
		\end{itemize}
	\end{mypro}
	\emph{Proof:}
	(i) can be obtained due to the conjugate theory, which is similar the proof of Proposition 1 in~\cite{WangZYRTSP2023}. As $\phi_{h,\lambda}$ is an even function and the function $g$ is convex, we then have (ii).
	For (iii), the bias $\Delta b$ can be estimated by the gap between the identity function and the proximity operator for $|x|\geq \lambda$, and since the proximity operator is odd, we only discuss $\Delta b= x-P_{\varphi_{h,\lambda}}(x)=a\cdot h'(x)$ for $x\geq \lambda$. Observing (\ref{Pro-L1}), it is apparent that the bias generated by the $\ell_1$-norm is $\lambda$. While when $x=\lambda$, the bias generated by our SIRs is $\lambda$, thus we need $\Delta b$ decreases with $x>\lambda$ to ensure that the bias produced by our SIRs is less than that by the $\ell_1$-norm. Therefore, we have $(\Delta b)'=a\cdot h''(x)<0$, implying that $\phi_{h,\lambda}(x)$ is strictly concave for $x>\lambda$. 
	$\hfill\blacksquare$
	
	\begin{table*}
		\caption{\small {Different loss functions and their proximity operators.}} 
		\vspace{-0.5cm}
		\begin{center}
			\setlength{\tabcolsep}{0.9mm}{
				{\begin{tabular}{ |c|c|c|c|}
						\hline
						& HOW & HOC & HOG\\
						\hline
						$l_{\cdot,\lambda}(x)$
						&$
						\begin{cases}
							x^2/2, &|x|\leq \lambda\\
							\frac{\sigma^2}{2}\left(1-e^{\frac{\lambda^2-x^2}{\sigma^2}}\right)+\frac{\lambda^2}{2}, &|x|\textgreater \lambda
						\end{cases}$
						&$
						\begin{cases}
							x^2/2, &|x|\leq \lambda\\
							\frac{\gamma^2+\lambda^2}{2}\ln \left(1+\frac{x^2}{\gamma^2}\right)+\delta, &|x|\textgreater \lambda
						\end{cases}$
						&$ 
						\begin{cases}
							x^2/2, &|x|\leq \lambda\\
							\frac{(\lambda^2+4\tau^2)^2x^2}{8\tau^2(x^2+4\tau^2)}-\frac{\lambda^4}{8\tau^2}, &|x|\textgreater \lambda
						\end{cases}$\\
						\hline
						$l_{\cdot,\lambda}(x)$
						&$\mathop {\min}\limits_{y} ~\frac{(x-y)^2}{2} + \lambda\varphi_{\sigma,\lambda}(y)$ 
						& $ \mathop {\min}\limits_{y} ~\frac{(x-y)^2}{2} + \lambda\varphi_{\gamma,\lambda}(y)$
						& $\mathop {\min}\limits_{y} ~\frac{(x-y)^2}{2} + \lambda\varphi_{\tau,\lambda}(y)$ \\
						\hline
						$P_{\varphi_{\cdot,\lambda}}(x)$ 
						& ${\rm max}\left\{0, |x|-|x|e^{\frac{\lambda^2-x^2}{\sigma^2}}\right\} {\rm sign}(x)$
						& ${\rm max}\left\{0, |x|-\frac{\left(\gamma^2+\lambda^2\right)|x|}{\gamma^2+x^2}\right\} {\rm sign}(x)$
						& ${\rm max}\left\{0, |x|-\frac{(\lambda^2+4\tau^2)^2|x|}{(x^2+4\tau^2)^2} \right\}{\rm sign}(x)$\\
						\hline
			\end{tabular}}}
			\label{regularization-funs}
		\end{center}
		\vspace{-0.5em}
	\end{table*}
	
\begin{figure}[htb]
		\centering
		\begin{minipage}{0.23\linewidth}
			\footnotesize
			\centerline{\includegraphics[width=4cm]{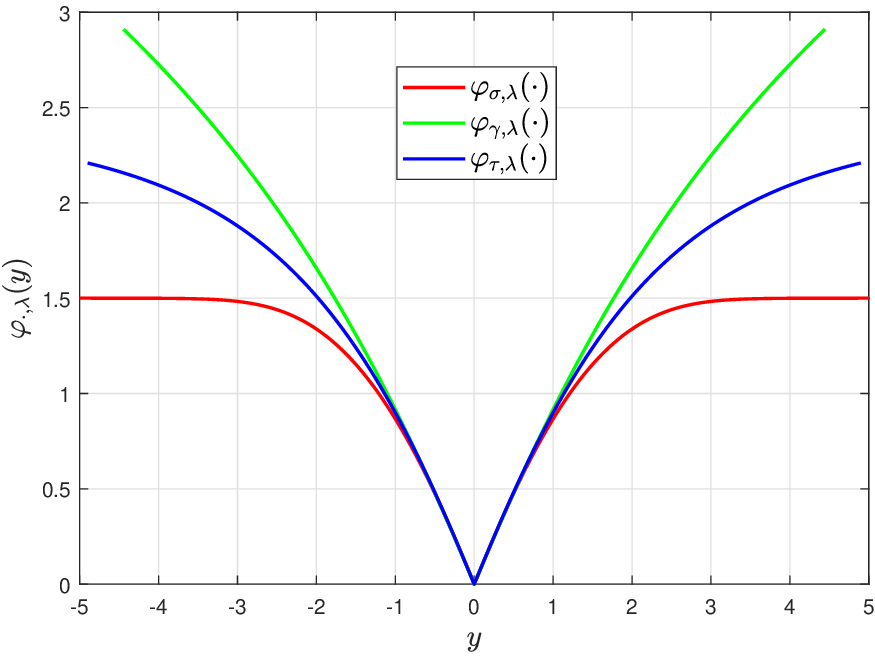}}\vskip 0pt
			\centerline{\scriptsize {(a)}}\vskip -3pt
		\end{minipage}\hspace{25mm}
		\begin{minipage}{0.24\linewidth}
			\footnotesize
			\centerline{\includegraphics[width=4cm]{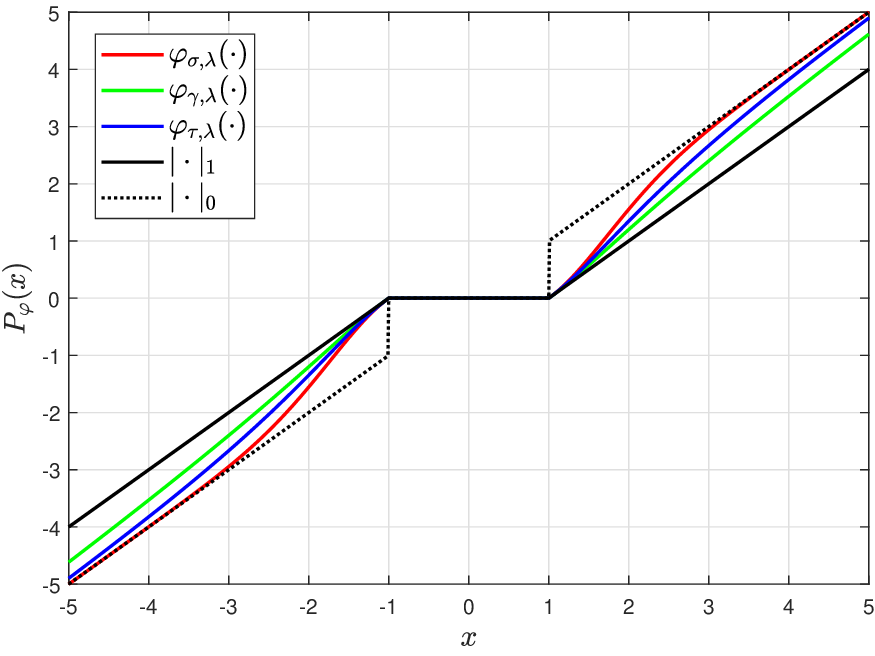}}\vskip 0pt
			\centerline{\scriptsize {(b)}}\vskip -3pt
		\end{minipage}\hspace{20mm}
		\vspace{-10pt}	
		\caption{Curves for (a) generated SIRs and (b) proximity operators with $\lambda=1$, $\sigma= \sqrt{2}\lambda$, $\gamma= \lambda$ and $\tau= {\sqrt{3}}\lambda/{2}$.}
		\label{proximity_curves}
	\end{figure}

	Next, we specify the function $g$ as the Welsch, GMC and Cauchy functions to develop the corresponding hybrid ordinary-Welsch (HOW), hybrid ordinary-GMC (HOG) and hybrid ordinary-Cauchy (HOC) functions, where  `ordinary' refers to the quadratic function, and the corresponding SIRs, namely, $\varphi_{\sigma,\lambda}(\cdot)$, $\varphi_{\gamma,\lambda}(\cdot)$ and $\varphi_{\tau,\lambda}(\cdot)$, which are shown in Table~\ref{regularization-funs}. Moreover, by Proposition~\ref{solution-proximal}, to make the bias generated by our SIRs less than that by the $\ell_1$-norm, we have $\sigma\leq \sqrt{2}\lambda$, $\gamma\leq \lambda$ and $\tau \leq \sqrt{3}\lambda/2$ for $\varphi_{\sigma,\lambda}(\cdot)$, $\varphi_{\gamma,\lambda}(\cdot)$ and $\varphi_{\tau,\lambda}(\cdot)$, respectively. 	
	Fig.~\ref{proximity_curves} plots the curves of $\varphi_{\sigma,\lambda}(\cdot)$, $\varphi_{\gamma,\lambda}(\cdot)$ and $\varphi_{\tau,\lambda}(\cdot)$ and their proximity operators. 
	
	\subsection{LRMC via Generated SIRs}
	
	In this section, we apply the generated SIRs to LRMC, resulting in: 
	\begin{equation}\label{Ori_pro}
		\begin{split}
			&\mathop {\min}\limits_{\pmb M}~ \|{\pmb M}\|_{\varphi_{\cdot,\lambda}}, ~ \text{s.t.} ~ \pmb M_{\Omega} = \pmb X_{\Omega}
		\end{split}
	\end{equation}
	where $\|\pmb M\|_{\varphi_{\cdot,\lambda}}=\sum_{i=1}^{r} \varphi_{\cdot,\lambda}{(\pmb \sigma_i)}$, and $\varphi_{\cdot,\lambda}(\cdot)$ can be the generated regularizer $\varphi_{\sigma,\lambda}(\cdot)$, $\varphi_{\gamma,\lambda}(\cdot)$ or $\varphi_{\tau,\lambda}(\cdot)$. 
	Problem (\ref{Ori_pro}) can be converted into the following equivalent form:
	\begin{equation}\label{RMC-formulation}
		\begin{split}
			&\mathop {\min}\limits_{\pmb M, \pmb E} ~\|\pmb M\|_{\varphi_{\cdot,\lambda}} ~\text{s.t.}~\pmb X = \pmb M + \pmb E,~ \pmb X_{\Omega^c}=0,~\pmb E_{ \Omega} = 0
		\end{split}
	\end{equation}
	which can be solved by the ADMM, and its augmented Lagrangian function is: 
	\begin{equation}\label{Aug_L}
		\begin{split}
			\mathcal{L}'_\rho(\pmb M, \pmb E,\pmb \Lambda) := &\|\pmb M\|_{\varphi_{\cdot,\lambda}} +\left<\pmb \Lambda, \pmb X-\pmb M -\pmb E\right>\\
			& +\frac{\rho}{2}\left\|\pmb X-\pmb M -\pmb E\right\|_F^2
		\end{split}
	\end{equation}
	which is equal to:
	\begin{equation}\label{Aug_L}
		\begin{split}
			\mathcal{L}_\rho(\pmb M, \pmb E,\pmb \Lambda) := &1/\rho\|\pmb M\|_{\varphi_{\cdot,\lambda}} +\frac{1}{2}\left\|\pmb X-\pmb M -\pmb E\right\|_F^2\\
			& +1/\rho\left<\pmb \Lambda, \pmb X-\pmb M -\pmb E\right>
		\end{split}
	\end{equation}
	where $\pmb \Lambda$ is the Lagrange multiplier matrix and $\rho>0$ is the penalty parameter. The variables $\pmb M, \pmb E$ and $\pmb \Lambda$ are updated at the $(k+1)$th iteration as follows:
	
	$Update~of $~$\pmb M$: Given $\pmb E^{k}$, $\pmb \Lambda^k$ and $\rho^k$, the estimated matrix $\pmb M$ is obtained by:  	
	\begin{equation}\label{update_M}
		{\rm arg}\mathop {\min}\limits_{\pmb M} 1/\rho^k\|\pmb M\|_{\varphi_{\cdot,1/\rho^k}} + \frac{1}{2}\left\|\pmb D^k -\pmb M\right\|_F^2, ~\pmb D^k = \pmb X-\pmb E^k + \frac{\pmb \Lambda^k}{\rho^k}
	\end{equation}
	If $\pmb D^k= \pmb U^k \Sigma^k \pmb V^k$ is the singular value decomposition (SVD) of $\pmb D^k$, then the optimal solution to (\ref{update_M}) according to Theorem 1 in~\cite{LuC2015} is:
	\begin{equation}\label{M_solution}
		\pmb M^{k+1}= \pmb U^k P_{\varphi_{\cdot,1/\rho^k}}(\pmb \Sigma^k) \pmb V^k
	\end{equation}
	
	$Update~of $~$\pmb E$: Given $\pmb M^{k+1}$, $\pmb \Lambda^k$ and $\rho^k$, $\pmb E^{k+1}$ is updated by solving:
	\begin{equation}\label{S_COmega}
		{\rm arg}\mathop {\min}\limits_{\pmb E_{\Omega^c}} \frac{1}{2}\left\|\pmb X_{\Omega^c}-\pmb M_{\Omega^c}^{k+1} + \frac{\pmb \Lambda_{\Omega^c}^k}{\rho^k} -\pmb E_{\Omega^c}\right\|_F^2
	\end{equation}
	with the optimal solution:
	\begin{equation}\label{S_solution_Om_S}
		\pmb E_{\Omega^c}^{k+1} = \frac{\pmb \Lambda_{\Omega^c}^k}{\rho^k} -\pmb M_{\Omega^c}^{k+1}
	\end{equation} 
	
	\begin{algorithm}[htb]
		\caption{MC via SIR}
		\label{Algo:NNSR}
		\algsetup{indent=1.5em}
		\vspace{1ex}
		\begin{algorithmic}
			\REQUIRE  $\pmb X_\Omega$, $\Omega$, $\mu=1.05$, $\xi=10^{-7}$ and $I_m=1000$
			\STATE \textbf{Initialize:} $\pmb E^0=\pmb 0$, $\pmb \Lambda^0=\pmb 0$, and $k=0$.
			\WHILE {$rel_E^k>\xi$ and $k\leq I_m$}
			
			\STATE Update $\pmb M^{k+1}$ via (\ref{M_solution})
			
			\STATE Update $\pmb E^{k+1}$ via (\ref{S_solution_Om_S})	
			
			\STATE Update $\pmb \Lambda^{k+1}= \pmb \Lambda^{k}+\rho^k\left(\pmb X -\pmb M^{k+1}-\pmb E^{k+1}\right)$
			
			\STATE Update $\rho^{k+1}=\mu\rho^k$	
			
			\STATE $k\leftarrow k+1$
			
			\ENDWHILE
			\ENSURE $\pmb M = \pmb M^{k}$.
		\end{algorithmic}
	\end{algorithm}
	
	The steps of the proposed algorithm are summarized in Algorithm~\ref{Algo:NNSR}. 	
	It is worth noting that the updates of $\pmb \Lambda$ and $\rho$ are only provided in Algorithm~\ref{Algo:NNSR} due to page limit. 
	When the generated regularizers $\varphi_{\sigma,\lambda}(\cdot)$, $\varphi_{\gamma,\lambda}(\cdot)$ are $\varphi_{\tau,\lambda}(\cdot)$ are adopted in (\ref{Ori_pro}), we refer the corresponding algorithm to as MC-HOW, MC-HOC and MC-HOG, respectively.
	We terminate Algorithm~1 until the relative error $ rel_E^k = \|\pmb X - \pmb M^k -\pmb S^k\|_F/\|\pmb X\|_F\leq \xi$ or the iteration number reaches the maximum allowable number $I_m$. 
	Besides, the SVD computation with complexity of $\mathcal{O}(\min (m,n)mn)$ is involved per iteration. Furthermore, the convergence analysis of Algorithm~\ref{Algo:NNSR} is provided in Proposition~\ref{Covergence-proof}, and its proof is similar to the proof of Theorem 3 in~\cite{GuSH2017}, thus we omit it due to page limit.
	\begin{mypro}\label{Covergence-proof} 
		The sequence $\{\pmb M^k, \pmb E^k, \pmb \Lambda^k\}$ generated in Algorithm 1 satisfies:
		\begin{itemize}
			\item[(i)] The generated sequences $\{\pmb M^k, \pmb E^k, \pmb \Lambda^k\}$ are bounded.
			\item[(ii)] $\lim_{k \to \infty}\left\|\pmb M^{k+1}-\pmb M^k\right\|_F^2=0$ \\
			$\lim_{k \to \infty}\left\|\pmb X-\pmb M^{k+1} -\pmb E^{k+1}\right\|_F^2=0$.
		\end{itemize}
	\end{mypro}

	\section{Experimental Results}
	\label{sec:sim}
	We compare our algorithms with the competing approaches, including NNM~\cite{LinZCMC2010}, IRNN-$\ell_p$ ($p=0.5$)~\cite{LuC2016} and IRNN-SCAD~\cite{LuC2014}. All numerical simulations are conducted using a computer with 3.0 GHz CPU and 16 GB memory. 
	\begin{figure}[htb]
		\centering
		
		\begin{minipage}{0.23\linewidth}
			\footnotesize
			\centerline{\includegraphics[width=4cm]{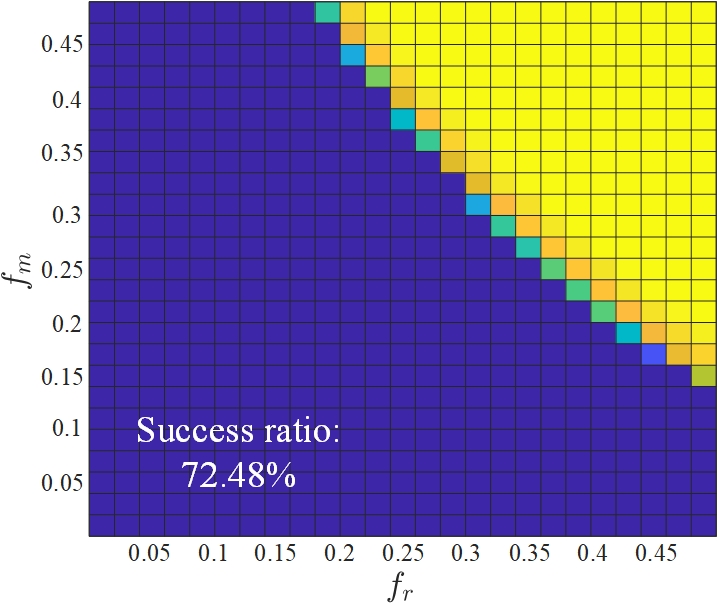}}\vskip 0pt
			\centerline{\scriptsize {NNM}}\vskip -3pt
		\end{minipage}\hspace{25mm}
		\begin{minipage}{0.24\linewidth}
			\footnotesize
			\centerline{\includegraphics[width=4.54cm]{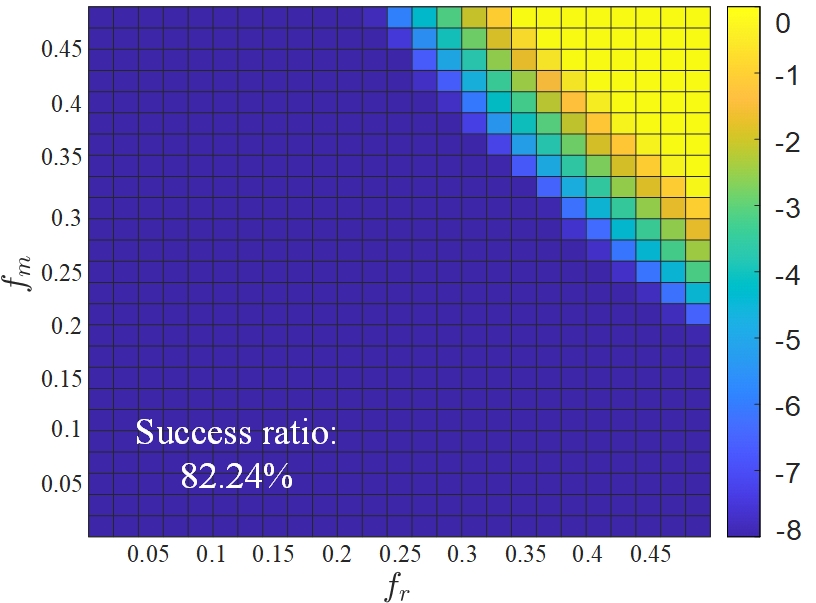}}\vskip 0pt
			\centerline{\scriptsize {IRNN-$\ell_p~(p=0.5)$}}\vskip -3pt
		\end{minipage}\hspace{20mm}
		
		\begin{minipage}{0.23\linewidth}
			\footnotesize
			\centerline{\includegraphics[width=4cm]{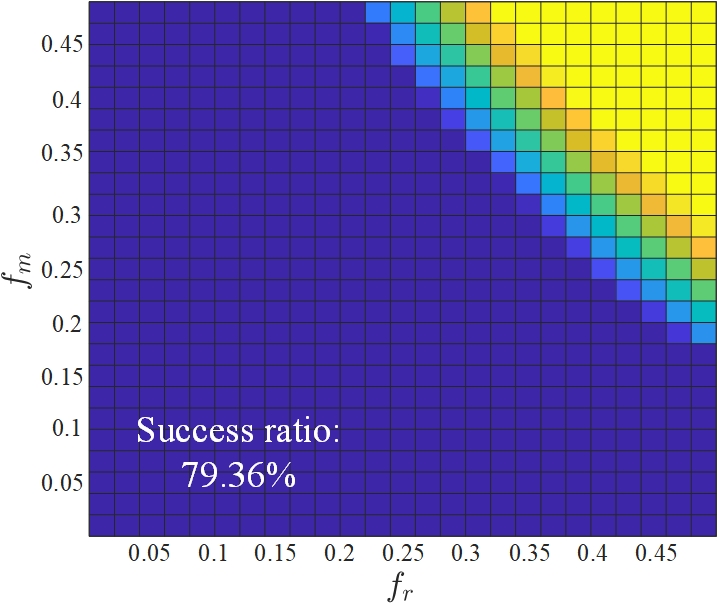}}\vskip 0pt
			\centerline{\scriptsize {IRNN-SCAD}}\vskip -3pt
		\end{minipage}\hspace{25mm}
		\begin{minipage}{0.24\linewidth}
			\footnotesize
			\centerline{\includegraphics[width=4.54cm]{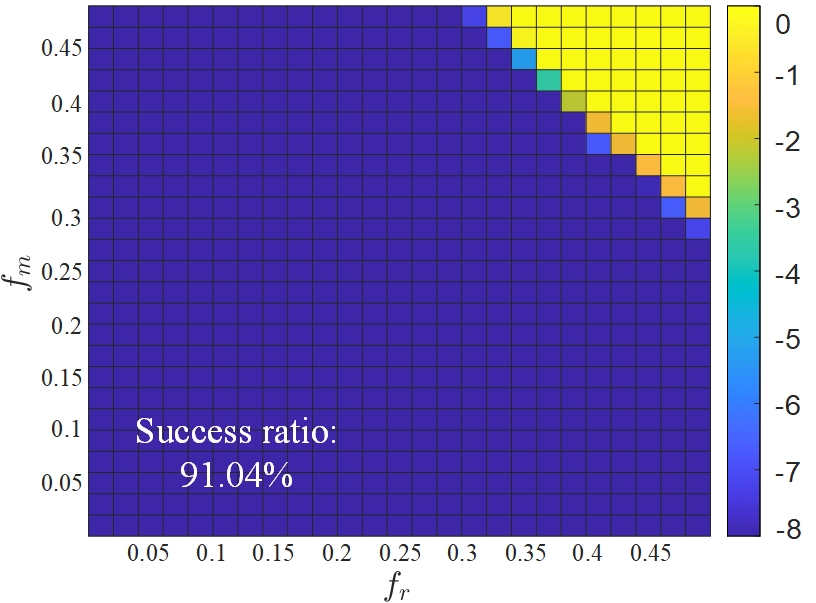}}\vskip 0pt
			\centerline{\scriptsize {MC-HOW}}\vskip -3pt
		\end{minipage}\hspace{20mm}
		
		\begin{minipage}{0.23\linewidth}
			\footnotesize
			\centerline{\includegraphics[width=4cm]{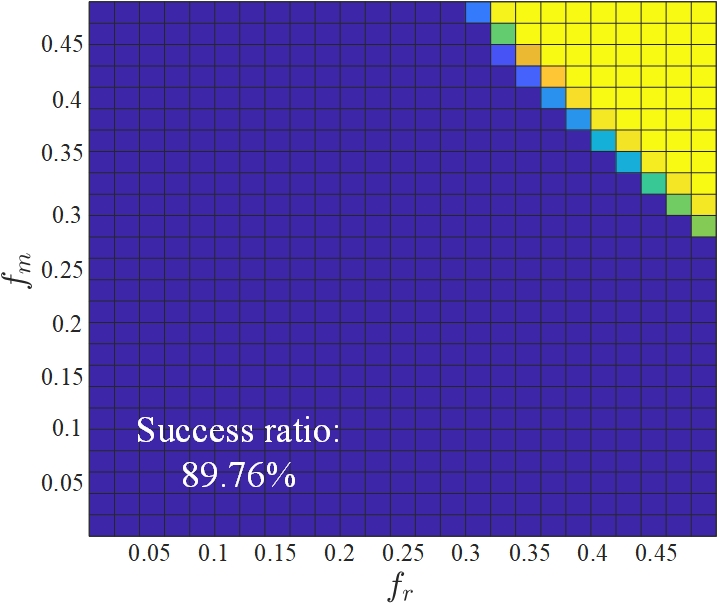}}\vskip 0pt
			\centerline{\scriptsize {MC-HOC}}\vskip -3pt
		\end{minipage}\hspace{25mm}
		\begin{minipage}{0.24\linewidth}
			\footnotesize
			\centerline{\includegraphics[width=4.54cm]{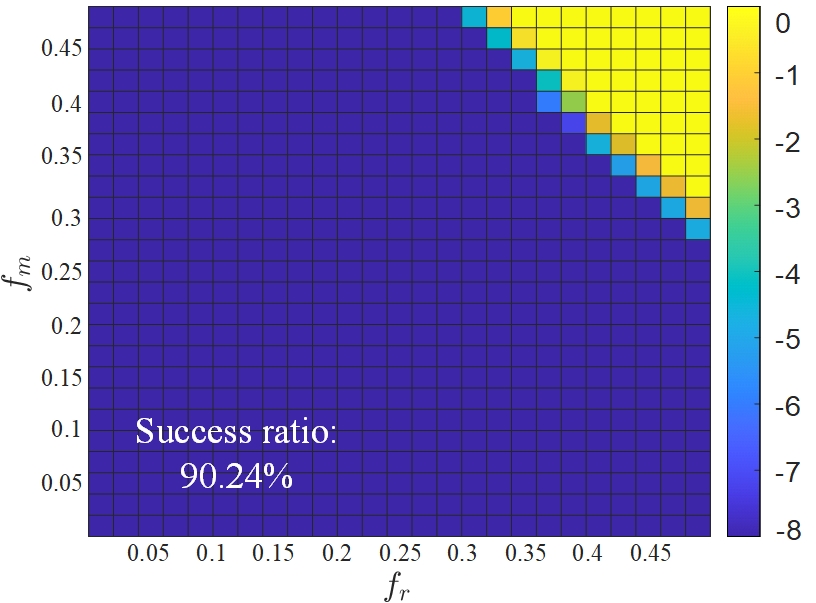}}\vskip 0pt
			\centerline{\scriptsize {MC-HOG}}\vskip -3pt
		\end{minipage}\hspace{20mm}
		
		\caption{Algorithm phase transition diagrams with different fractions of full rank and missing entries $\{f_r,f_m\}$. }
		\label{log_RRE_map}
	\end{figure}
	A rank-$r$ matrix $\pmb X = \pmb U \pmb V$ is first generated, where the entries of $\pmb U \in \mathbb{R}^{m\times r}$ and $\pmb V \in \mathbb{R}^{r\times n}$ are sampled from the standard Gaussian distribution. In this study, $r=f_r \times n$ where $f_r$ is the fraction of full rank. Besides, we randomly remove $f_m\times m\times n$ entries from $\pmb X$, where $f_m$ is the fraction of missing entries, to 
	yield the incomplete matrix $\pmb X_\Omega$. Furthermore, root mean square error (RMSE) defined as ${\rm RMSE} = {\|\pmb X - \pmb M\|_F}/{\sqrt{mn}}$ is adopted to measure the recovery performance.
	In our experiments, $m=300$, $n=200$, $f_m$ and $f_r$ range from $0.01$ to $0.05$ with a step size of $0.02$. All methods are evaluated by the RMSE based on $100$ independent runs.
	Fig.~\ref{log_RRE_map} shows the log-scale RMSE with different $\{p_r,p_s\}$ values.
	It is seen that compared with the convex NNM, the nonconvex rank surrogates can recover more cases, and MC-HOW yields the biggest success area. Here, if RMSE $<10^{-3}$, we consider that is a success recovery. 
	
	\begin{figure}[htb]
		\centering
		\includegraphics[width=6cm]{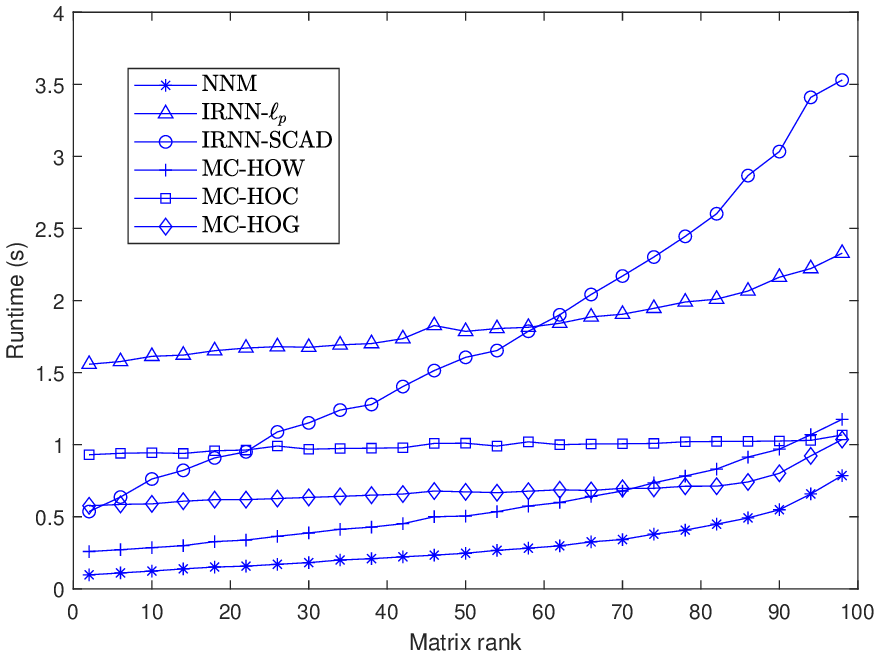}
		\vspace{-0.5em}
		\caption{Runtime versus matrix rank.}\label{Runtime_competitors}
	\end{figure}
	
	On the other hand, the runtime (in seconds) for all techniques is investigated under different matrix ranks. Fig.~\ref{Runtime_competitors} plots runtime versus matrix rank with $f_m=0.1$. We see that the runtime of our approaches is less than that of IRNN-$\ell_p$ because our regularizers have closed-form proximity operators. Nevertheless, NNM is the most computationally efficient since the proximity operator of the $\ell_1$-norm has a simpler expression than those of our SIRs.
	
	\section{Conclusion}
	\label{sec:con}
	
	In this paper, we devise a framework to generate SIRs with closed-form proximity operators. We analyze that the Moreau envelope of our regularizers is a convex problem although the regularizers may be nonconvex, and provide the corresponding closed-form solution. Besides, it is proved that the bias generated by our SIRs is less than that by the $\ell_1$-norm under certain conditions. Then, we employ our SIRs as nonconvex rank surrogates for LRMC, and algorithms based on the ADMM are developed.  
	Finally, extensive numerical experiments are conducted to demonstrate that the developed algorithms can achieve better recovery under different matrix ranks and missing ratios. 
		
	\bibliographystyle{IEEEbib}
	\bibliography{strings,refs}

\end{document}